\documentclass[11pt,leqno,twoside]{article}

\usepackage{amsfonts,amsmath,amsthm,amssymb}
\usepackage{enumerate}
\usepackage{graphics,graphicx,subfigure}

\usepackage{color}
\usepackage{todonotes}
\usepackage{cancel}
\usepackage{url}
\usepackage{hyperref}
\usepackage{makeidx}
\usepackage{showidx}
\usepackage{multicol}        
\usepackage{xspace}
\usepackage{stmaryrd}        
\usepackage{pifont}          
\usepackage{fancybox}        
\usepackage{bm}

\usepackage{fancyhdr}
\theoremstyle{plain}
\theoremstyle{definition}
\newtheorem{lem}{Lemma}
\newtheorem{defn}[lem]{Definition}
\newtheorem{thm}[lem]{Theorem}
\newtheorem{prop}[lem]{Proposition}
\newtheorem{cor}[lem]{Corollary}
\newtheorem{notn}[lem]{Notations}
\newtheorem{pb}[lem]{Problem}
\newtheorem{form}[lem]{Formulae}

\newtheorem*{rk}{Remark}
\newtheorem*{com}{Comment}
\newtheorem*{ex}{Example}
\theoremstyle{remark}

\newcommand{\blem}{\begin{lem}}
\newcommand{\elem}{\end{lem}}
\newcommand{\bdefn}{\begin{defn}}
\newcommand{\edefn}{\end{defn}}
\newcommand{\bthm}{\begin{thm} }
\newcommand{\ethm}{\end{thm}}
\newcommand{\bprop}{\begin{prop}}
\newcommand{\eprop}{\end{prop}}
\newcommand{\bcor}{\begin{cor}}
\newcommand{\ecor}{\end{cor}}
\newcommand{\bnotn}{\begin{notn}}
\newcommand{\enotn}{\end{notn}}
\newcommand{\bpb}{\begin{pb}}
\newcommand{\epb}{\end{pb}}
\newcommand{\bform}{\begin{form}}
\newcommand{\eform}{\end{form}}
\newcommand{\brk}{\begin{rk}}
\newcommand{\erk}{\end{rk}}
\newcommand{\bcom}{\begin{com}}
\newcommand{\ecom}{\end{com}}
\newcommand{\bex}{\begin{ex}}
\newcommand{\eex}{\end{ex}}
\newcommand{\bpf}{\begin{proof}}
\newcommand{\epf}{\end{proof}}
\newcommand{\rev}[1]{{{#1}}} 

\newcommand{\vQ}{{\bf Q}}

\newcommand{\cA}{\mathcal{A}}

\newcommand{\cH}{\mathcal{H}}

\newcommand{\cK}{\mathcal{K}}
\newcommand{\cL}{\mathcal{L}}

\newcommand{\cO}{\mathcal{O}}
\newcommand{\cP}{\mathcal{P}}

\newcommand{\cU}{\mathcal{U}}
\newcommand{\cV}{\mathcal{V}}
\newcommand{\cW}{\mathcal{W}}
\newcommand{\cX}{\mathcal{X}}
\newcommand{\cY}{\mathcal{Y}}

\newcommand{\bC}{\mathbb{C}}
\newcommand{\bD}{\mathbb{D}}

\newcommand{\bN}{\mathbb{N}}

\newcommand{\bR}{\mathbb{R}}

\newcommand{\bT}{\mathbb{T}}

\newcommand{\be}{\begin{equation}}
\newcommand{\ee}{\end{equation}}
\newcommand{\bal}{\begin{align}}
\newcommand{\eal}{\end{align}}
\newcommand{\ba}{\begin{align*}}
\newcommand{\ea}{\end{align*}}
\newcommand{\bmx}{\begin{matrix}}
\newcommand{\emx}{\end{matrix}}
\newcommand{\bbmx}{\begin{bmatrix}}
\newcommand{\ebmx}{\end{bmatrix}}
\newcommand{\bpmx}{\begin{pmatrix}}
\newcommand{\epmx}{\end{pmatrix}}
\newcommand{\bvmx}{\begin{vmatrix}}
\newcommand{\evmx}{\end{vmatrix}}

\newcommand{\ol}{\overline}

\newcommand{\wt}{\widetilde}
\newcommand{\f}{\frac}

\newcommand{\inc}{\subseteq}

\newcommand{\setm}{\setminus}

\newcommand{\argmin}{{\rm argmin}\,}

\newcommand{\minimize}[1]{\underset{#1}{\rm minimize}\,}

\renewcommand{\Re}{\operatorname{Re}}

\newcommand{\la}{\lambda}

\newcommand{\eps}{\varepsilon}

\title{\vspace{-20mm}Approximability Models and Optimal System Identification\medskip\hrule height 1.2pt \vspace{-6mm}}
\author{Mahmood Ettehad and Simon Foucart\footnote{S. F. partially supported by NSF grants DMS-1622134 and DMS-1664803.} --- Texas A\&M University}
\date{\vspace{-6mm}\rule{100mm}{0.8pt}}

\newcommand\shorttitle{Approximability Models and Optimal System Identification}
\newcommand\authors{M. Ettehad, S. Foucart}

\begin{document}
\maketitle

\vspace{-0mm}
\begin{abstract}
This article considers the problem of optimally recovering stable linear time-invariant \mbox{systems} observed via linear measurements made on their transfer functions.
A common modeling \mbox{assumption} is replaced here by the related assumption that the transfer functions belong to a model set described by approximation capabilities. Capitalizing on recent optimal-recovery results relative to such approximability models, we construct some optimal algorithms and \mbox{characterize} the optimal performance for the identification and evaluation of transfer functions in the framework of the Hardy Hilbert space and of the disc algebra.
In particular, we determine explicitly the optimal recovery performance for frequency measurements taken at equispaced points on an inner circle or on the torus. 
\end{abstract}

\noindent {\it Key words and phrases:}  optimal recovery, system identification, Hardy spaces, disc algebra.

\noindent {\it AMS classification:} 93B30, 46J10, 30H10.

\vspace{-5mm}
\begin{center}
\rule{100mm}{0.8pt}
\end{center}


\section{Background and Motivation}

System identification can be viewed as the learning task of inferring a model to describe the behavior of a system observed via input-output data. As such, it influences control theory as a whole, 
since the model selection dictates the subsequent design of efficient controllers for the system. Linear dynamical systems, due to their appearance in many applications,
have been the focus of many investigations.
They are often handled by moving
to the frequency domain, where one studies their transfer functions,
typically assumed to be elements of the Hardy space $\cH_2$ (see e.g. \cite[Chapter 13]{Zhou}) or $\cH_\infty$ (see e.g. \cite[Chapter 14]{Zhou}).
Then, based on {\em a priori} information given as an hypothesized model and on {\em a posteriori} information given as frequency observations, one seeks to identify the transfer function in a way that minimizes the error with respect to the given norm.

For the model put forward in \cite{HJN},
various algorithms for the identification of transfer functions in $\cH_\infty$ have already been proposed:
in \cite{CNF92}, interpolatory algorithms are constructed by solving a Nevanlinna--Pick problem, \cite{CNFan92} presents some closely related algorithms with explicit bounds in $\cH_\infty$ and $\ell_1$ norms, and in \cite{CN93}  interpolatory algorithms are generalized to time-domain data by solving a Carath\rev{\'e}odory--Fej\rev{\'e}r problem via standard convex optimization methods.
More recent works focus on the trade-off between the number of samples required to accurately build models of dynamical systems and the degradation of performance in various control objectives.
For instance, the article \cite{TBPR} derived bounds on the number of noisy input-output samples from a stable linear time-invariant system that are sufficient to guarantee 
closeness of the finite impulse response approximation to the true system in $\mathcal{H}_{\infty}$-norm. 
\rev{In \cite{SBTR}, linear system identification from pointwise noisy measurements in the frequency domain was formulated as a convex minimization problem involving an $\ell_1$-penalty and error estimates in $\cH_2$-norm were provided.}
There has also been progress in quantifying the sample complexity of dynamical system identification using \rev{statistical estimation theory,} see \cite{VRK08,CW02}.

\vspace{-2mm}
The article \cite{HJN} was one of the first to suggest adopting a perspective from optimal recovery \cite{MicRiv} in system identification, leading to a framework that is now textbook material (see e.g. \cite[Section 4.4]{Part} or \cite[Section~3.2]{ChenGu}).
The purpose of this note is to revisit this classical framework in light of recent optimal-recovery results involving approximability models, see \cite{BCDDPW,DPW,DFPW}.
Our setting is closely related to the one from \cite{HJN}, in that
the {\em a priori} information is encapsulated by a model set $\cK$
and the {\em a posteriori} information consists of frequency response measurements (possibly nonuniformly spaced, as in \cite{A94}), while the objective is to establish performance bounds and devise algorithms for the recovery of transfer functions $F$.
By recovering $F$, we mean here approximating it in full (i.e., identifying) in the context of $\cH_2$
or evaluating a point value $F(\zeta_0)$ (i.e., estimating) in the context of $\cH_\infty$ ---  in fact, as in \cite{BSG95}, we replace $\cH_\infty$ by a subspace known as the disc algebra to avoid certain technical difficulties.
The novelty of our setting lies in the model set $\cK$:
whereas the model set of \cite{HJN} can be viewed as an intersection $\cK = \cap_{n \ge 0} \cK_n$ of approximation sets,
we focus here on a single approximation set $\cK_n$ chosen as model set.
This slight modification of the framework allows us to construct identification algorithms that are {\em optimal},
not just near-optimal.
In addition, these optimal algorithms turn out to also be {\em linear} algorithms.

\vspace{-2mm}
The structure of this article is as follows.
In Section \ref{SecPF}, we formulate precisely the problem we are considering, we introduce the approximability model,
and we describe some optimal-recovery results that have recently emerged. In Section \ref{SecH2}, we zoom in on results about Hilbert spaces and adapt them to the identification of transfer functions in $\cH_2$, leading to a matrix-analytic method for the construction of an optimal algorithm and the determination of the optimal performance. In the case of data gathered at equispaced points on inner circles, we also find the exact value of the optimal performance in terms of the number $m$ of observations and we remark that it is independent of the dimension $n$  of the polynomial space underlying the approximability model.
In Section \ref{SecDA}, we turn to results about quantities of interest in Banach spaces and adapt them to the optimal evaluation of transfer functions in the disc algebra, leading to a convex-optimization method for the construction of an optimal algorithm and the determination of the optimal performance.
In the case of data gathered at equispaced points on the torus,
we also show how the optimal performance behaves in terms of $m$ and $n$. Section \ref{SecConc} concludes with some perspective on further research.
Finally, an appendix collects some that have been delayed to keep the flow of the main text going.

\section{Problem Formulation}
\label{SecPF}

In its most abstract form,
the scenario we are considering in the rest of this article involves unknown objects $F$ from a normed space $\cX$ which are acquired through observational data
\be
y_k = \ell_k(F), \qquad k= 1,\dots,m,
\ee
where $\ell_1,\ldots,\ell_m \in \cX^*$ are known linear functionals.
The perfect accuracy of the acquisition process is of course an idealization. For brevity, we shall write $y = \cL(F) \in \bC^m$.
Discarding some of the $\ell_k$'s if necessary,
we may assume that the operator $\cL:\cX \to \bC^m$ has full range.
The task at hand consists in making use of the data $y$ to approximate $F$, or merely to evaluate a quantity of interest $Q(F)$,  where $Q$ is a linear map from $\cX$ into another normed space $\cY$ --- typically $\cY = \bC$.
There is also an {\em a priori} knowledge about $F$,
often conveyed by the assumption that $F$  belongs to a certain model set $\cK \inc \cX$. The feasibility of the task is assessed in a worst-case setting via the quantity
\be
\label{Opt}
E^{\rm opt}(\cK,\cL,Q) := \inf_{A: \bC^m \to \cY}  \sup_{F \in \cK}
\| Q(F) - A(\cL(F)) \|_\cY.
\ee
An optimal algorithm (relative to the model set $\cK$) is a map $A^{\rm opt}$ from $\bC^m$ into $\cY$ for which the infimum is achieved.

In system identification, the standard objects are transfer functions
belonging to the Hardy space $\cX = \cH_2(\bD)$ or $\cX=\cH_\infty(\bD)$ relative to the open unit disk $\bD := \{ z \in \bC: |z|<1 \}$.
We recall that the Hardy spaces $\cH_p(\bD)$ are defined for $1 \le p \le \infty$ by
\be
\cH_p(\bD) := \left\{ F \mbox{ is analytic on }\bD \mbox{ and } \|F\|_{\cH_p} < + \infty \right\},
\ee
where the $\cH_p$-norms of $F(z) = \sum_{n=0}^\infty f_n z^n$ are given for $p=2$ and $p=\infty$ by
\begin{eqnarray}
\|F\|_{\cH_2}  := 
& \displaystyle{\sup_{r \in [0,1)}} \bigg[ \displaystyle{\f{1}{2 \pi} \int_0^{2 \pi} |F(r e^{i \theta})|^2 d \theta} \bigg]^{1/2}
& = \bigg[\sum_{n=0}^\infty |f_n|^2 \bigg]^{1/2},\\
\label{DefHInfNorm}
\|F\|_{\cH_\infty} := 
& \displaystyle{\sup_{r \in [0,1)}} \quad \sup_{\theta \in [-\pi, \pi]} |F(r e^{i \theta})| \quad \qquad
& = \sup_{ |z|=1 } |F(z)|.
\end{eqnarray}
The last equality in \eqref{DefHInfNorm} does not imply that functions in $\cH_\infty(\bD)$ are well defined on the torus $\bT := \{ z \in \bC: |z| = 1 \}$. To bypass this peculiarity,
instead of the whole Hardy space $\cH_\infty(\bD)$,
we shall work within the subspace consisting of functions that are continuous on $\bT$. This space is called the disc algebra and is denoted $\cA(\bD)$.

Concerning the acquisition process, although our considerations are valid for any linear functionals $\ell_1,\ldots,\ell_m$,
we shall concentrate in this initial work on the situation commonly encountered in system identification where $\ell_1,\ldots,\ell_m$ are point evaluations at some $\zeta_1,\ldots,\zeta_m$ located  in the disc or on the torus.
It is worth recalling that the point evaluation defined at some $\zeta \in \bC$ by
\be
e_\zeta(F) := F(\zeta),
\qquad F \in \cX,
\ee
is well defined for $|\zeta| \le 1$ when $\cX = \cA(\bD)$,
but only for $|\zeta| < 1$ when $\cX = \cH_2(\bD)$.
Its Riesz representer $E_\zeta \in \cH_2(\bD)$, 
characterized by the identity $e_\zeta(F) = \langle F,E_\zeta \rangle_{\cH_2}$ for all $F \in \cH_2(\bD)$, is the Cauchy kernel given by
\be
\label{RepPtEv}
E_\zeta(z) := \f{1}{1-\ol{\zeta} z},
\qquad |z|<1.
\ee
We also point out that, in $\cA(\bD)$, the dual norm of a point evaluation at $\zeta \in \bC$, $|\zeta| \le 1$, is
\be 
\label{NormPE}
\|e_{\zeta}\|_{\cA^*} = 1,
\ee
where $\|e_{\zeta}\|_{\cA^*} \ge 1$ follows e.g. from $e_\zeta(1) = 1$ and where $\|e_{\zeta}\|_{\cA^*} \le 1$ holds because,
for any $F \in \cA(\bD)$, one has $|F(\zeta)| \le \|F\|_{\cH_\infty}$ by the maximum modulus theorem.

As for the model, a popular representation of the {\em a priori} information proposed in \cite{HJN} prevails.
Given $\rho > 1$ and $M>0$, it in fact involves two models sets, $\cK_{\cH_2}$ and $\cK_{\cH_\infty}$,
both of them consisting of functions $F$ that are analytic on the disc  $\bD_\rho := \{ z \in \bC: |z|< \rho \}$ and that satisfy, respectively,
\begin{align}
\label{Mod2}
\|F(\rho \, \cdot)\|_{\cH_2} & = \bigg[\sum_{n=0}^\infty |f_n|^2 \rho^{2n} \bigg]^{1/2} \le M,\\
\label{ModInf}
\|F(\rho \, \cdot)\|_{\cH_\infty} &= \sup_{|z| = \rho} |F(z)| \le M. 
\end{align}
\rev{These model sets are closely related
	(see the appendix for the precise statement and its justification)
	to the model sets given,
	for $\cX = \cH_2(\bD)$ or $\cX = \cH_\infty(\bD)$, by
}
\be
\label{KasIntersect}
\rev{ \bigcap_{n \ge 0}
	\left\{ F \in \cX: {\rm dist}_{\cX}(F,\cP_n) \le \eps_n \right\},}
\ee
where $\cP_n$ denotes the space of polynomials of degree at most $n-1$ and where $\eps_n = M \rho^{-n}$ 
(with constants $\rho>1$ and $M>0$ differing from the ones in \eqref{Mod2}-\eqref{ModInf}). In the rest of this article,
we consider the model sets obtained by overlooking the intersection and focusing on single sets associated with fixed $n$'s.
In this situation, we will be able to produce optimal algorithms in the sense of \eqref{Opt}.
As we will realize, these optimal algorithms are in fact linear maps\footnote{\rev{In the case $\cK = \cH_2(\bD)$, the optimal algorithm over an ellipsoidal model set such as the one described by \eqref{Mod2} is linear, too. It is a variant of the minimal-norm interpolation presented in Section 5.3. of \cite{Part}}}.
Our method leverages recent results described \rev{next}.
Strictly speaking, they were established in \cite{BCDDPW,DPW,DFPW} for the real setting only --- their validity for the complex setting is justified   in the appendix.

Inspired by parametric PDEs, the authors of the articles \cite{BCDDPW,DPW,DFPW} emphasized approximation properties and considered the model sets
\be
\cK_\cX(\eps,\cV) := \left\{
F \in \cX: {\rm dist}_\cX(F,\cV) \le \eps
\right\}
\ee
relative to a normed space $\cX$, a subspace $\cV$ of $\cX$,
and a  parameter $\eps > 0$. As a summary  of general results for the full approximation problem where $Q={\rm Id}$ (see \cite{DPW} for details), we underline that
\vspace{-5mm}
\begin{enumerate}[(i)]
	\item \label{i} an important role is played by an indicator of the compatibility of the acquisition process and the model 
	(as represented by the spaces $\ker(\cL)$ and $\cV$, respectively),
	namely by
	\be
	\label{mu}
	\mu_\cX(\ker(\cL),\cV)
	:= \sup_{F \in \ker(\cL)} \f{\|F\|_{\cX}}{{\rm dist}_{\cX}(F,\cV)};
	\ee
	\item \label{ii} the optimal performance over the model set $\cK_\cX(\eps,\cV)$ is essentially characterized by this quantity, since
	\be
	\label{Ineqs}
	\mu_\cX(\ker(\cL), \cV) \; \eps 
	\le
	E^{\rm opt}(\cK_\cX(\eps,\cV),\cL,I) 
	\le 
	2 \, \mu_\cX(\ker(\cL), \cV) \; \eps ;
	\ee
	\item \label{iii} the map $A': \bC^m \to \cX$ defined independently of $\eps>0$ by
	\be
	\label{GenAlgo}
	A'(y) := \underset{F \in \cX}{\rm argmin} \; {\rm dist}_\cX(F,\cV)
	\qquad \mbox{subject to } \cL(F) = y
	\ee
	provides a near-optimal algorithm (though not necessarily a practical one) in the sense that
	\be
	\sup_{F \in \cK_\cX(\eps,\cV)} \|F-A'(\cL(F))\|_\cX \le 2 \, \mu_\cX(\ker(\cL),\cV) \; \eps .
	\ee
\end{enumerate}
It is worth noting that $\mu_\cX(\ker(\cL),\cV)$ decreases when observations are added, as $\ker(\cL)$ shrinks, and increases when the space $\cV$ is enlarged, as ${\rm dist}_\cX(F,\cV)$ becomes smaller. \rev{In fact, $\mu_\cX(\ker(\cL),\cV) = \infty$ when $n:= \dim(\cV) > m$ since then one can find $F \in \ker (\cL) \cap \cV$.}
Thus, we assume that $n \le m$ throughout. The articles \cite{BCDDPW} and \cite{DFPW} improved the results \eqref{i}-\eqref{ii}-\eqref{iii} in two specific situations.
Precisely, \cite{BCDDPW} considered the full approximation problem when $\cX$ is a Hilbert space, while \cite{DFPW} dealt with an arbitrary normed space $\cX$ but placed the focus on  quantities of interest $Q$ that are linear functionals. 
Our objective consists in adapting and supplementing these contributions for the spaces $\cX = \cH_2(\bD)$ and $\cX = \cA(\bD)$,
which is done in Section \ref{SecH2} and Section \ref{SecDA}, respectively.
\rev{More precisely, our novel contribution consists,
	for the case of $\cH_2(\bD)$,
	in a slight but computation-ready variation on the result of \cite{BCDDPW}
	(Theorem \ref{ThmH2})
	and in the exact determination of the indicator $\mu_{\cH_2}(\ker(\cL_{\zeta}),\cP_n)$
	when $\zeta_1,\ldots,\zeta_m$ are equispaced points on an inner circle (Proposition \ref{PropMuH2Equi}),
	and for the case of $\cA(\bD)$,
	in a slight but appropriate extension of the result of \cite{DFPW} 
	(Theorem \ref{ThmHInf})
	and in the asymptotic determination of 
	the indicator $\mu_{\cA}(\ker(\cL_{\zeta}),\cP_n)$
	when $\zeta_1,\ldots,\zeta_m$ are equispaced points on the torus (Proposition \ref{PropMuHinfEqui}).
}

\section{Optimal Identification in $\cH_2(\bD)$}
\label{SecH2}

As was just mentioned, the results \eqref{i}-\eqref{ii}-\eqref{iii} can be enhanced when $\cX$ is a Hilbert space, which is the case of the Hardy space $\cH_2(\bD)$.
To present the enhancement provided by \cite{BCDDPW},
we let $(V_1,\ldots,V_n)$ denote an orthonormal basis for the subspace $\cV$ of the Hilbert space $\cX$
and $L_1,\ldots,L_m \in \cH_2(\bD)$ denote the Riesz representers of  the linear functionals $\ell_1,\ldots,\ell_m \in \cX^*$ --- recall that they are characterized by
\be
\ell_k(F) = \langle F, L_k \rangle_{\cH_2},
\qquad k=1,\ldots,m.
\ee
Concerning \eqref{iii},
it was shown that the algorithm of \eqref{GenAlgo} is not just near optimal, but genuinely {\em optimal}. 
It can in fact be written as the {\em linear} map $A^{\rm opt}: \bC^m \to \cX$ given by
\be
\label{OptAlgH2}
A^{\rm opt}(y) = V^\star_y+ W_y- P_\cW(V^\star_y),
\qquad  \qquad
V^\star_y= \underset{V \in \cV}{\argmin} \|W_y-P_\cW(V)\|_{\cX},
\ee
where $P_\cW$ is the orthogonal projector onto the space $\cW := {\rm span}\{ L_1,\ldots,L_m \}$ and $W_y$  denotes the element of $\cW$ satisfying $\ell_k(W_y) = y_k$ for all $k=1,\ldots,m$ 
(so that $W_y= P_\cW(F)$ if $y= \cL(F)$ for some $F \in \cX$).
Concerning \eqref{ii}, the optimal performance over $\cK_{\cX}(\eps,\cV)$ is then completely characterized by a sharpening of \eqref{Ineqs} to 
\be
\label{OptPerfH2}
E^{\rm opt}(\cK_\cX(\eps,\cV),\cL,I) = \mu_\cX(\ker(\cL),\cV) \; \eps.
\ee
Concerning \eqref{i}, the compatibility indicator introduced in \eqref{mu} becomes fully computable as
\be
\label{MuH2}
\mu_\cX(\ker(\cL),\cV) = \f{1}{\sigma_{\min}(\wt{G})},
\ee
where $\sigma_{\min}(\wt{G})$ is the smallest positive singular value of the cross-Gramian matrix $\wt{G} \in \bC^{m \times n}$ relative to the orthonormal basis $(V_1,\ldots,V_n)$ for $\cV$
and to an orthonormal basis $(\wt{L}_1,\ldots,\wt{L}_m)$~for~$\cW$. 
The entries of this matrix are 
\be
\wt{G}_{k,j} := \langle  V_j, \wt{L}_k \rangle,
\qquad k=1,\ldots,m, \quad j=1,\ldots,n.
\ee

\subsection{Optimal performance and algorithm}

When turning to the implementation of these results for the Hardy space $\cH_2(\bD)$, a difficulty  arises from the fact that orthonormal bases for $\cW$ are not automatically available.
A Gram--Schmidt orthonormalization of $(L_1,\ldots,L_m)$ is not ideal, because calculating inner products in $\cH_2$ is not exact 
(if performed either as contour integrals or as inner products of infinite sequences).
We are going to present a more reliable method of implementation,
based only on matrix computations and relying on data directly available to the users.
These data consist of two matrices $G \in \bC^{m \times n}$ and $H \in \bC^{m \times m}$.
The first one, with entries
\be
\label{DefG}
G_{k,j} := \ell_k(V_j),
\qquad k=1,\ldots,m, \quad j=1,\ldots,n,
\ee
is simply the cross-Gramian matrix relative to $(V_1,\ldots,V_n)$ and to $(L_1,\ldots,L_m)$, 
in view of the identity $\ell_k(V_j) = \langle V_j, L_k \rangle_{\cH_2}$.
The second one, with entries
\be
\label{DefH}
H_{k,j} := \ell_k(L_j),
\qquad k,j = 1,\ldots,m,
\ee
is the Gramian matrix relative to $(L_1,\ldots,L_m)$,
in view of the identity $\ell_k(L_j) = \langle L_j, L_k \rangle_{\cH_2}$.

Once the matrices $G$ and $H$ are set,
we can make very explicit the optimal performance of system identification in $\cH_2$ for the approximability model relative to a subspace $\cV$, as well as an algorithm achieving this optimal performance.
\rev{This is the object of the theorem below --- to reiterate, it is a variation on a result of \cite{BCDDPW} which advantageously lends itself more easily to practical implementation.}

\vspace{3mm}
\bthm
\label{ThmH2}
With $G \in \bC^{m \times n}$ and $H \in \bC^{m \times m}$ defined in \eqref{DefG} and \eqref{DefH},
one has
\be
\label{muH2}
\mu_{\cH_2}(\ker(\cL),\cV) =  \f{1}{\la_{\min}(G^* H^{-1} G)^{1/2}}.
\ee
Moreover,
with $c^\star = (G^* H^{-1} G)^{-1}  G^* H^{-1} y$
and $d^\star = H^{-1}(y- G c^\star)$,
the map $A^{\rm opt}: \bC^m \to \cH_2(\bD)$ defined by 
\be
\label{AoptH2}
A^{\rm opt}(y) = 
\sum_{j=1}^n c^\star_j V_j + \sum_{k=1}^m d^\star_k L_k
\ee
is an optimal algorithm in the sense that
\be
\sup_{F \in \cK_{\cH_2}(\eps,\cV)} \|F - A^{\rm opt}(\cL(F)) \|_{\cH_2}
= \inf_{A: \bC^m \to \cH_2} 
\sup_{F \in \cK_{\cH_2}(\eps,\cV)} \|F - A(\cL(F)) \|_{\cH_2},
\ee
with the value of the latter being $\mu_{\cH_2}(\ker(\cL),\cV)  \, \eps$.
\ethm

\bpf Since the Gramian matrix $H$ is positive definite,
it has an eigenvalue decomposition 
\be
\label{eigH}
H = U D U^*,
\ee
where $U \in \bC^{m \times m}$ is a unitary matrix and $D = {\rm diag}[d_1,\ldots,d_m] \in \bC^{m \times m}$ is a diagonal matrix with positive entries.
It is routine to verify that 
the functions $\wt{L}_1,\ldots,\wt{L}_m \in \cH_2(\bD)$ defined by
\be 
\label{DefLTilde}
\wt{L}_k = \f{1}{\sqrt{d_k}} \sum_{j=1}^m U_{j,k} L_j,
\qquad
k=1,\ldots,m,
\ee 
form an orthonormal basis for $\cW$. It is also easy to verify that the cross-Gramian matrix $\wt{G} \in \bC^{m \times n}$ relative to $(V_1,\ldots,V_n)$ and to $(\wt{L}_1,\ldots,\wt{L}_m)$
is then expressed as  
\be
\label{Gt}
\wt{G} = D^{-1/2} U^* G.
\ee
In turn, we derive that
\be
\wt{G}^* \wt{G} = G^* H^{-1} G.
\ee
The first part of the theorem, namely \eqref{muH2},
is now a consequence of \eqref{MuH2}. For the second part of the theorem, by virtue of \eqref{OptAlgH2}, our strategy is simply to determine in turn $W_y$, $V^\star_y$, and $P_\cW (V^\star_y)$. To start with, it is readily checked, taking inner products with $L_1,\ldots,L_k$, that 
\be
W_y= \sum_{k=1}^m (H^{-1}y)_k L_k. 
\ee
With the othonormal basis $(\wt{L}_1,\ldots,\wt{L}_m)$ for $\cW$ introduced in \eqref{DefLTilde},
the easily verifiable change-of-basis formula
\be
\label{CofB}
W = \sum_{k=1}^m c_k L_k 
\iff
W = \sum_{k=1}^m (D^{1/2} U^* c)_k \wt{L}_k
\ee
yields the alternative representation
\be 
\label{CofB2}
W_y= \sum_{k=1}^m (D^{-1/2} U^* y)_k \wt{L}_k.
\ee
Next, for an arbitrary function $V = \sum_{j=1}^n c_j V_j \in \cV$,
we observe that
\be
\label{PW}
P_\cW(V) = \sum_{k=1}^m \langle V, \wt{L}_k \rangle \wt{L}_k
= \sum_{k=1}^m \sum_{j=1}^n c_j \langle V_j, \wt{L}_k\rangle \wt{L}_k
= \sum_{k=1}^m (\wt{G} c)_k\wt{L}_k. 
\ee
It follows from \eqref{CofB2} and \eqref{PW} that
\be
\|W_y- P_\cW(V)\|_{\cH_2}^2
= \sum_{k=1}^m |(D^{-1/2} U^* y)_k - (\wt{G} c)_k |^2
= \| D^{-1/2} U^* y- \wt{G} c \|_2^2.
\ee
Therefore, the minimizer $V^\star_y\in \cV$ of this expression takes the form
\be
\label{Vstar}
V^\star_y= \sum_{j=1}^n c^\star_j V_j,
\qquad
c^\star = \wt{G}^\dagger D^{-1/2} U^* y,
\ee
where $\wt{G}^\dagger = (\wt{G}^* \wt{G})^{-1} \wt{G}^*$ is the Moore--Penrose pseudo-inverse of $\wt{G}$.
According to \eqref{Gt}, we easily derive that the coefficient vector $c^\star$ is indeed
\be
\label{Cstar}
c^\star = (G^* H^{-1} G)^{-1}  G^* H^{-1} y.
\ee
Finally, in view of \eqref{PW} and \eqref{CofB}, we obtain
\be
\label{PVstar}
P_\cW(V^\star_y) = \sum_{k=1}^m ( \wt{G} c^\star )_k \wt{L}_k
= \sum_{k=1}^m (U D^{-1/2} \wt{G} c^\star)_k L_k
= \sum_{k=1}^m (H^{-1} G c^\star)_k L_k.
\ee
Putting \eqref{Vstar}-\eqref{Cstar}, \eqref{CofB2}, and \eqref{PVstar} together in \eqref{OptAlgH2},
we arrive at the expression announced in \eqref{AoptH2}.
This completes the proof of the theorem.
\epf

\subsection{Evaluations at points on an inner circle}

It is time to specify the results to our scenario of interest.
In particular, keeping \eqref{KasIntersect} in mind,
we fix $\cV$ to be the space $\cP_n$ of polynomials of degree at most $n-1$. It is equipped with the orthonormal basis $(V_1,\ldots,V_n)$ given by $V_j(z)= z^{j-1}$, $j=1,\ldots,n$.
Let us also consider a radius $r<1$ and suppose that the linear functionals $\ell_k$, $k=1,\ldots,m$, take the form 
\be
\ell_k(F) = F(\zeta_k)
\qquad \mbox{for some }\zeta_k \in \bC \mbox{ with } |\zeta_k| = r. 
\ee
We are going to compare two situations,
one where $\zeta_1,\ldots,\zeta_m$ are randomly selected on the circle of radius $r$ and one where they are equispaced on this circle. In the latter situation, the optimal performance of system identification in $\cH_2$ can be determined explicitly.

\bprop
\label{PropMuH2Equi}
Given $0 < r < 1$ and $\zeta_1,\ldots,\zeta_m \in \bD$ defined by 
\be
\zeta_k = r \exp \left( i \f{2\pi }{m} (k-1) \right),
\qquad k=1,\ldots,m,
\ee
the indicator of compatibility of the acquisition process $\cL_{\zeta}$ and the approximability model relative to $\cP_n$ is,
for any $n \le m$,
\be
\label{MuEquiH2}
\mu_{\cH_2}(\ker(\cL_{\zeta}),\cP_n) = \f{1}{\sqrt{1-r^{2m}}}.
\ee
\eprop 

\bpf
We shall prove the stronger statement that
\be
\label{DiagForm}
G^* H^{-1} G = (1-r^{2m}) I_n,
\ee
which clearly implies the announced result by taking \eqref{muH2} into account. To this end, we notice that the entries of the matrix $G \in \bC^{m \times n}$ are
\be
\label{ExprG}
G_{k,j} = V_j(\zeta_k) = r^j \exp \left( i \f{2 \pi}{m}(k-1)(j-1) \right),
\qquad k=1,\ldots,m, 
\quad j = 1,\ldots,n,
\ee
and, in view of the form \eqref{RepPtEv} of the representer $L_j = E_{\zeta_j}$,
that the entries of the matrix $H \in \bC^{m \times m}$ are
\be
H_{k,j} = L_j(\zeta_k) = \f{1}{1-\ol{\zeta_j} \zeta_k} = \f{1}{1-r^2 \exp \left( i \f{2 \pi (k-j)}{m} \right)},
\qquad k=1,\ldots,m, 
\quad
j = 1,\ldots,m.
\ee
Since $H$ is a circulant matrix, it `diagonalizes in Fourier',
meaning that the matrix $U$ in the eigendecomposition \eqref{eigH} has columns $u^{(1)},\ldots,u^{(m)} \in \bC^m$ given by
\be
u^{(k)} = \f{1}{\sqrt{m}}
\bbmx
1\\
\exp(i 2 \pi (k-1)/m)\\
\vdots \\
\exp(i 2 \pi (k-1)(m-1)/m)
\ebmx .
\ee
The eigenvalues $d_1,\ldots,d_m$ are found through the calculation 
\begin{align}
(H u^{(k)})_{k'} 
& = \sum_{j=1}^{m} \f{1}{1-r^2 \exp( i 2 \pi (k'-j)/m )} \f{\exp(i 2 \pi  (k-1)(j-1) )}{\sqrt{m}}\\
\nonumber
& = \sum_{j=1}^{m} \f{\exp(i 2 \pi  (k-1)(j-k')  / m)}{1-r^2 \exp( i 2 \pi (k'-j)/m )} (u^{(k)})_{k'},
\end{align}
so that, after the change of summation index $h=k'-j$, 
the eigenvalue $d_k$ appears to be
\begin{align}
d_k & = \sum_{h=1}^{m} \f{\exp(-i 2 \pi (k-1) h / m)}{1-r^2 \exp( i 2 \pi h/m )}
= \sum_{h=1}^{m} \exp(-i 2 \pi  (k-1)h / m) \sum_{t=0}^\infty r^{2t} \exp( i 2 \pi  t h/m )
\\ \nonumber
& = \sum_{t=0}^\infty r^{2t} \sum_{h=1}^{m} \exp( i 2 \pi  (t-k+1)h/m )
= \sum_{t=0}^\infty r^{2t} \, m \, {\bf 1}_{\{t \in k-1+ m \bN \}} 
\\ \nonumber
& = m (r^{2(k-1)} + r^{2(k-1+m)} + r^{2(k-1+2m)} + \cdots)
= \f{m r^{2(k-1)}}{1-r^{2m}}.
\end{align}
Thus, we have now made explicit the eigendecomposition of $H$ as
\be
\label{FormH}
H = \f{m}{1-r^{2m}} \, U \,  {\rm diag} [1,r^2,\ldots,r^{2(m-1)}] \, U^*,
\qquad \quad U = \bbmx u^{(1)}  |  \cdots  |  u^{(m)} \ebmx.
\ee
Besides, we easily observe that the expression \eqref{ExprG} reads, in matrix form,
\be
\label{FormG}
G = \sqrt{m}\,  \wt{U} \, {\rm diag} [1,r,\ldots,r^{n-1}],
\qquad \quad \wt{U} = \bbmx u^{(1)} | \cdots | u^{(n)} \ebmx. 
\ee
At this point, it is an easy matter to verify  that \eqref{FormH} and \eqref{FormG} imply \eqref{DiagForm}. 
\epf

\begin{figure}[htbp]
	\center
	\subfigure[The indicator $\mu_{\cH_2}(\ker(\cL_{\zeta}),\cP_n)$ as a function of $n$ for equispaced points (no dependence) and for random points (fast increase).]{
		\includegraphics[width=0.45\textwidth]{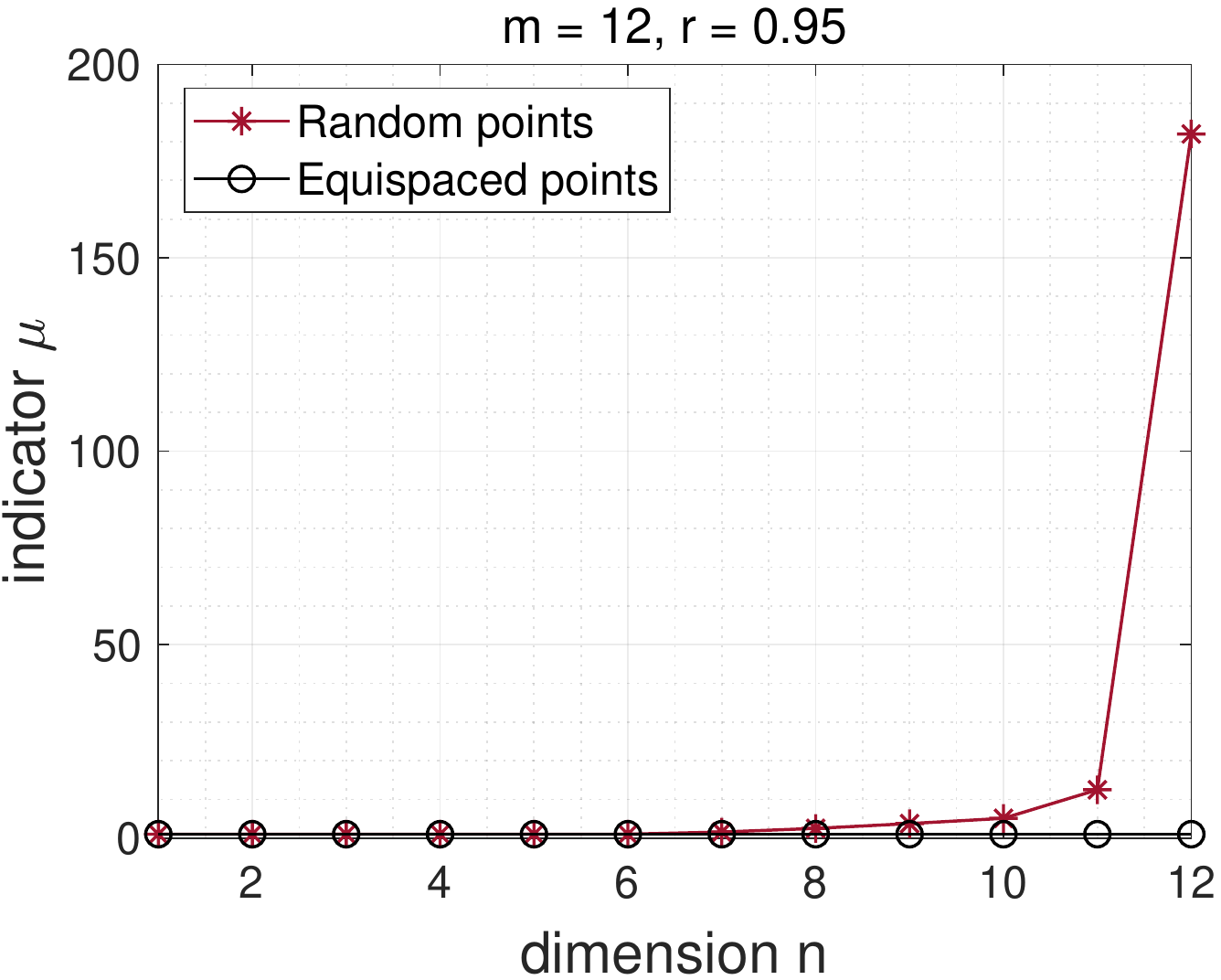}
	}
	\hspace{5mm}
	\subfigure[The identification error $\|F-A^{\rm opt}(\cL_\zeta(F))\|_{\cH_2}$ as a function of $n$ for equispaced points and for random points.]{
		\includegraphics[width=0.45\textwidth]{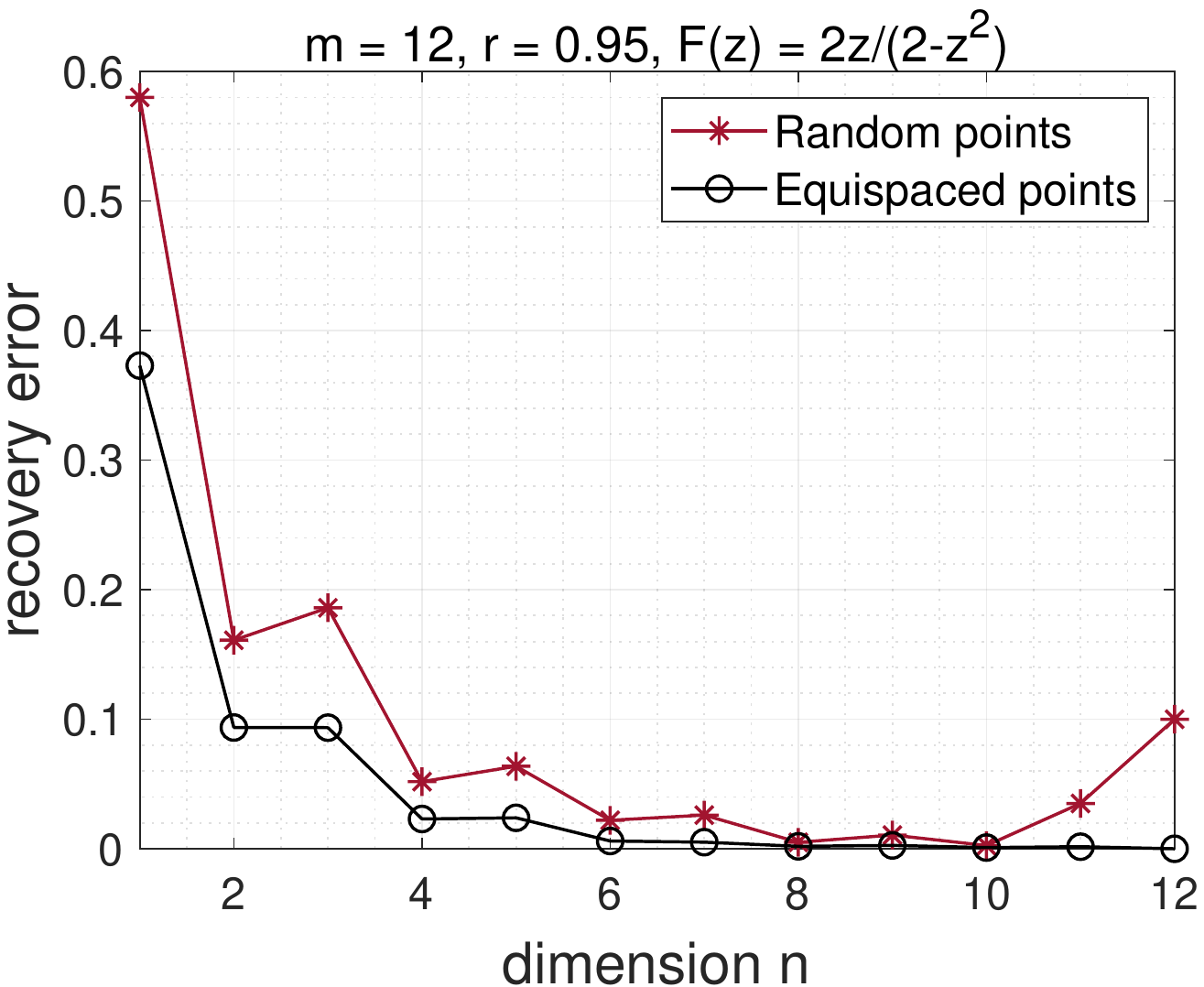}
	}
	\caption{Behavior of the indicator $\mu_{\cH_2}(\ker(\cL_{\zeta}),\cP_n)$ and of the identification error \mbox{$\|F - A^{\rm opt}(\cL_\zeta(F))\|_{\cH_2}$} when
		the points $\zeta_1,\ldots,\zeta_m$ are chosen on a circle of radius $r<1$.}
	\label{Fig1}
\end{figure}

In the case of observations gathered at equispaced points,
Proposition \ref{PropMuH2Equi} shows that the indicator $\mu_{\cH_2}(\ker(\cL_{\zeta}),\cP_n)$ is independent of the dimension $n$ of $\cV = \cP_n$. Keeping in mind that the value of $\eps=\eps_n$ appearing in \eqref{OptPerfH2} decreases with $n$
(e.g. as $\eps_n = M \rho^{-n}$ for the model set imposed by \eqref{Mod2}-\eqref{KasIntersect}), the optimal performance becomes most favorable when $n$ as large as possible, i.e., when $n=m$.
This differs from the typical situation where a tradeoff is to be found between the increase of $\mu_{\cH_2}(\ker(\cL_{\zeta}),\cP_n)$ and the decrease of $\eps_n$. 
Figure \ref{Fig1} illustrates numerically the fact that $\mu_{\cH_2}(\ker(\cL_{\zeta}),\cP_n)$ generally increases with $n$.
The experiments that generated this figure can be reproduced by downloading from the authors' webpages the {\sc matlab} files accompanying this article.\footnote{To compute the $\cH_2$-error between functions $F = \sum_{j=1}^\infty b_j V_j$ and $\wt{F} = \sum_{j=1}^n c_j V_j + \sum_{k=1}^m d_k L_k$,
	we used the fact that $\| F -\wt{F} \|_{\cH_2}^2 = \|F\|_{\cH_2}^2 + \|\wt{F}\|_{\cH_2}^2 - 2 \Re \langle F, \wt{F} \rangle $,
	together with $\|\wt{F}\|_{\cH_2}^2 = \|c\|_2^2 + \langle d, H d \rangle + 2 \Re \langle c, G d \rangle$ and $\langle F, \wt{F} \rangle = \langle b_{1:n}, c \rangle + \sum_{k=1}^m  \ol{d_k} \ell_k(F)$.}

\brk
\rev{The choice $n=m$ implies the
	existence of some $A_m : \bC^m \to \cH_2(\bD)$ yielding the estimate
	\be
	\label{NearBest}
	\|F - A_m(\cL_\zeta(F))\|_{\cH_2}
	\le \f{1}{\sqrt{1-r^{2m}}} \, {\rm dist}_{\cH_2}(F,\cP_m)
	\qquad 
	\mbox{for all } F \in \cH_2(\bD),
	\ee
	where the factor $1/\sqrt{1-r^{2m}}$ is bounded independently of $m$ by $1/\sqrt{1-r^{2}}$.}
As it turns out, the map $F \mapsto A_m(\cL_\zeta(F))$ is just the operator $\rev{P_m}: \cH_2(\bD) \to \cP_m$ of polynomial interpolation at $\zeta_1,\ldots,\zeta_m$:
indeed, thanks to $\cL_\zeta(\rev{P_m}(F)) = \cL_\zeta(F)$,
\eqref{NearBest} applied to $\rev{P_m}(F) \in \cP_m$ instead of $F$ reveals that $A_m(\cL_\zeta(F)) = \rev{P_m}(F)$.
In short, the operator $\rev{P_m}: \cH_2(\bD) \to \cP_m$ of polynomial interpolation at $m$ equispaced points on the circle of radius \rev{$r <1$}
acts as an operator of near-best approximation from $\cP_m$ relative to $\cH_2$.
\rev{As pointed out by a referee,
	this specific observation can be obtained without relying on the preceding  machinery.
	Indeed,
	the expression
	\be
	P_m(G)(z) = \sum_{\ell=0}^{m-1} \Big( \sum_{t=0}^\infty g_{\ell+tm} r^{tm} \Big) z^{\ell}
	\qquad
	\mbox{whenever }
	G(z) = \sum_{\ell=0}^\infty g_{\ell} z^{\ell}  
	\ee
	is easily verified by checking that $P_m(G)(\zeta_k) = G(\zeta_k)$ for all $k=1,\ldots,m$
	and that $P_m(V_j) = V_j$ for all $j=1,\ldots,m$.
	Then an application of Cauchy--Schwarz inequality gives
	\be
	\|P_m(G)\|_{\mathcal{\cH}_2} \le \f{1}{\sqrt{1-r^{2m}}} \|G\|_{\mathcal{\cH}_2}.
	\ee
	The inequality \eqref{NearBest} finally follows by taking $G = F-V$
	where $V \in \cP_m$ is chosen such that $\|F-V\|_{\mathcal{H}_2} = {\rm dist}_{\cH_2}(F,\cP_m)$.
	We emphasize that this short argument for \eqref{NearBest} incidentally justifies that 
	$\mu_{\cH_2}(\ker(\cL_{\zeta}),\cP_m) \le 1/\sqrt{1-r^{2m}}$,
	according to the leftmost inequality in  \eqref{Ineqs},
	and in turn that $\mu_{\cH_2}(\ker(\cL_{\zeta}),\cP_n) \le 1/\sqrt{1-r^{2m}}$ for all $n \le m$.
}
\erk

\section{Optimal Estimation and Identification in $\cA(\bD)$}
\label{SecDA}

As mentioned at the end of Section \ref{SecPF},
the results \eqref{i}-\eqref{ii}-\eqref{iii} can also be enhanced when $\cX$ is an arbitrary normed space and one does not target the full approximation of $F \in \cX$ but only the estimation of a quantity of interest $Q(F)$ for some linear functional $Q \in \cX^*$.
To present the enhancement provided in  \cite{DFPW},
we let $(V_1,\ldots,V_n)$ denote a basis for the subspace $\cV$ of $\cX$ --- this time, it does not need to be an orthonormal basis.
Concerning \eqref{iii}, an algorithm which is {\em optimal} over the approximation set $\cK_\cX(\eps,\cV)$ was uncovered:
it consists in outputting
\be
\label{OptAlgHInf}
A^{\rm opt}(y) = \sum_{k=1}^m a^\star_k y_k,
\qquad 
y\in \bC^m,
\ee
where the vector $a^\star \in \bC^m$ is obtained as a solution of the optimization problem
\be 
\label{OptinX*}
\minimize{a \in \bC^m} \bigg\| Q - \sum_{k=1}^m a_k \ell_k \bigg\|_{\cX^*}
\qquad \mbox{subject to }
\sum_{k=1}^m a_k \ell_k(V_j) = Q(V_j) \quad \mbox{for all }j=1,\ldots, n.
\ee
Note that the vector $a^\star$ is computed offline once and for all and is subsequently used for each fresh occurrence of observational data $y$ via the rule \eqref{OptAlgHInf}, thus providing an optimal algorithm which is incidentally a {\em linear} functional.
Note also that the knowledge of $\eps$ is unnecessary to produce the vector $a^\star$. Concerning \eqref{ii}, or rather its alteration to incorporate quantities of interest, the optimal performance over $\cK_\cX(\eps,\cV)$ is again completely characterized by
\be
E^{\rm opt}(\cK_\cX(\eps,\cV),\cL,Q) = \mu_\cX(\ker(\cL),\cV,Q) \, \eps,
\ee 
where the indicator of compatibility of the acquisition process and the model now becomes
\be
\label{IndicAD}
\mu_\cX(\ker(\cL),\cV,Q) :=
\sup_{F \in \ker(\cL)} \f{|Q(F)|}{{\rm dist}_\cX(F,\cV)}.
\ee
Concerning \eqref{i},
the value of this indicator is in fact obtained as the minimum of the optimization program \eqref{OptinX*}. At first sight, this program seems intractable because the dual norm is not automatically amenable to numerical computations. The purpose of this section is to show that the optimization can be performed effectively when the acquisition functionals are evaluations at points on the torus $\bT$ and when $\cX$ is the disc algebra $\cA(\bD)$ endowed with the $\cH_\infty$-norm.
It could also be performed when $\cX = \cH_2(\bD)$
--- we do not give details about this case, since it is essentially covered in \cite[Subsection 5.2]{DFPW}, which is dedicated to reproducing kernel Hilbert spaces and hence applies to $\cH_2(\bD)$.

\subsection{Optimal performance and algorithm}

The key to transforming \eqref{OptinX*} into a practically solvable optimization program when $\cX = \cA(\bD)$ consists in reformulating the $\cX^*$-norm of a linear combination of linear functionals as a workable expression of its coefficients.
With acquisition functionals being evaluations at points on the torus, this reformulation is made possible by Rudin--Carleson theorem. It provides access not only to the optimal performance of system identification in $\cA(\bD)$ for the approximability model relative to a subspace $\cV$, but also to an algorithm achieving this optimal performance.
\rev{This is the object of the theorem below
	--- to reiterate, it extends a result of \cite{DFPW} from the space of continuous functions to the more involved disc algebra.}

\bthm
\label{ThmHInf}
Given distinct points $\zeta_1,\ldots,\zeta_m \in \bT$
and a quantity of interest taking the form
$Q(F) =  F(\zeta_0)$ for some $\zeta_0 \in \bT \setm \{\zeta_1,\ldots,\zeta_m\}$,
consider the $\ell_1$-minimization problem
\be
\label{OptinAD}
\minimize{a \in \bC^m} \sum_{k=1}^m |a_k|
\qquad \mbox{subject to }
\sum_{k=1}^m a_k V_j(\zeta_k) = V_j(\zeta_0) \quad \mbox{for all }j=1,\ldots, n.
\ee 
If $a^\star \in \bC^m$ denotes a solution to this problem, then
\be
\mu_{\cA} (\ker(\cL_{\zeta}), \cV, Q)
= 1+\sum_{k=1}^m |a^\star_k|,
\ee
and the linear functional $A^{\rm opt}: \bC^m \to \bC$ defined by
\be
A^{\rm opt}(y) = \sum_{k=1}^m a^\star_k y_k,
\qquad y \in \bC^m,
\ee
is an optimal algorithm over $\cK_{\cA}(\eps,\cV)$ in the sense that
\be
\sup_{F \in \cK_{\cA}(\eps,\cV)} \rev{|} Q(F) - A^{\rm opt}(\cL_{\zeta}(F)) \rev{|} 
= \inf_{A: \bC^m \to \bC} 
\sup_{F \in \cK_{\cA}(\eps,\cV)} \rev{|} Q(F) - A(\cL_{\zeta}(F)) \rev{|},
\ee
with the value of the latter being $\mu_{\cA} (\ker(\cL_{\zeta}), \cV, Q)  \, \eps.$
\ethm

\bpf
The argument is based on the remark that, for any $a \in \bC^m$,
\be
\label{NormHInf*asL1}
\bigg\| e_{\zeta_0}  - \sum_{k=1}^m a_k e_{\zeta_k} \bigg\|_{\cA^*}
= 1+ \sum_{k=1}^m |a_k|.
\ee
The `$\le$'-part follows from the triangle inequality and the observation \eqref{NormPE}. As for the `$\ge$'-part, since the set $\{\zeta_0,\zeta_1,\ldots, \zeta_m \} \inc \bT$
is closed and have measure zero, Rudin--Carleson theorem (see e.g. \cite[Theorem~2.3.2]{Part}) ensures that one can find $F \in \cA(\bD)$ such that 
\be
\label{RC}
\|F\|_{\cH_\infty} = 1,
\qquad
F(\zeta_0) = 1 =: w_0,
\qquad
F(\zeta_k) =  - \f{\ol{a_k}}{|a_k|} =: w_k,
\quad k = 1,\ldots, m,
\ee
which in turn implies that 
\be
\bigg\| e_{\zeta_0}  - \sum_{k=1}^m a_k e_{\zeta_k} \bigg\|_{\cA^*}
\ge   F(\zeta_0)  - \sum_{k=1}^m a_k F(\zeta_k)
= 1 + \sum_{k=1}^m |a_k|.
\ee
It is then clear from \eqref{NormHInf*asL1} that the optimization program \eqref{OptinX*} reformulates as \eqref{OptinAD},
omitting the additive constant $1$. The announced results are now restatements of \eqref{OptAlgHInf}-\eqref{OptinX*}
and of the announced equality between the indicator \eqref{IndicAD} and the minimum of the optimization program.
\epf

\brk
The argument \eqref{RC} justifying the identity \eqref{NormHInf*asL1} would not hold if $\zeta_0,\zeta_1,\ldots, \zeta_m$ were all inside the unit circle.
Indeed, by Nevanlinna--Pick theorem (see e.g. \cite[Theorem 2.1.6]{Part} or \cite[Theorem 2.3.4]{ChenGu}),
it would mean that the matrix $\vQ \in \bC^{(m+1)\times (m+1)}$ with entries
\be 
Q_{k,j} = \f{1-\ol{w_j} w_k}{1-\ol{\zeta_j} \zeta_k},
\qquad k,j = 0,1,\ldots,m,
\ee
is positive semidefinite.
This cannot occur because its diagonal entries are all zero, hence its trace equals zero, 
so $\vQ$ itself would be the zero matrix.
\erk

\brk
Theorem \ref{ThmHInf} can effortlessly be extended to quantities of interest taking the form of a weighted sum of evaluations at points on $\bT$.
\rev{With a little more work,
	it could also be extended to quantities of interest of the form $Q(F) = F^{(s)}(\zeta)$ for some $\zeta \in \bD$.
	The details are left to the reader.
	They essentially amount to justifying the `$\ge$'-part of an analog of \eqref{NormHInf*asL1}.
	This is done via a limiting argument which involves a discretization of the Cauchy formula for $F^{(s)}(\zeta)$ and allows for Rudin--Carleson theorem to be applied. 
}
\erk

\subsection{Evaluations at points on the unit circle}

In this subsection, general considerations are again particularized to our situation of interest where we fix $\cV$ to be the space of polynomials of degree at most $n-1$. In this case, and with acquisition functionals being evaluations at equispaced points on the torus $\bT$, we are able to determine quite explicitly the optimal performance of system identification in $\cA(\bD)$.
Underlying the argument is a crucial observation about $\mu_{\cA}(\ker(\cL_{\zeta}),\cV)$ valid regardless of the space $\cV$ and of the evaluation points $\zeta_1,\ldots,\zeta_m \in \bT$.
This observation is isolated below.

\blem
\label{LemObs}
Given a finite-dimensional subspace $\cV$ of $\cA(\bD)$
and distinct points $\zeta_1,\ldots,\zeta_m \in \bT$, 
\be
\label{mu=1+}
\mu_{\cA}(\ker(\cL_{\zeta}),\cV)
= 1 + \sup_{V \in \cV} \f{\|V\|_{\cH_\infty}}{\max\limits_{k=1,\ldots,m} |V(\zeta_k)|}.
\ee
\elem

\bpf
In order to establish \eqref{mu=1+}, we first recall that 
\be
\mu_{\cA}(\ker(\cL_{\zeta}),\cV)
= \sup_{F(\zeta_1)=\cdots=F(\zeta_m)=0} \f{\|F\|_{\cH_\infty}}{{\rm dist}_{\cH_\infty}(F,\cV)}.
\ee
To prove the `$\le$'-part of \eqref{mu=1+},
we remark that, if $F \in \cA(\bD)$ satisfies $F(\zeta_1)=\cdots = F(\zeta_m) = 0$ and $V \in \cV$ satisfies $\|F - V \|_{\cH_\infty} = {\rm dist}_{\cH_\infty}(F,\cV)$,
then
\begin{align}
\f{\|F\|_{\cH_\infty}}{{\rm dist}_{\cH_\infty}(F,\cV)}
& = \f{\|F\|_{\cH_\infty}}{\|F - V \|_{\cH_\infty} }
\le 1 + \f{\|V\|_{\cH_\infty}}{\|F - V \|_{\cH_\infty} }
\le 1 + \f{\|V\|_{\cH_\infty}}{\max\limits_{k=1,\ldots,m} |(F-V)(\zeta_k)|}\\
\nonumber
& = 1 + \f{\|V\|_{\cH_\infty}}{\max\limits_{k=1,\ldots,m} |V(\zeta_k)|},
\end{align}
and the required inequality immediately follows.
As for the `$\ge$'-part, let $\eta$ denote the maximum appearing in \eqref{mu=1+},
and let us select $V \in \cV$ with $\max\limits_{k=1,\ldots,m}|V(\zeta_k)| = 1$
and $\|V\|_{\cH_\infty} = \eta$.
We then pick $z \in \bT$ such that $|V(z)| = \eta$
--- we may assume $z \not\in \{ \zeta_1,\ldots,\zeta_m \}$,
otherwise we can slightly perturb it  and replace $\eta$ by $\eta - \eps$ for an arbitrary small $\eps >0$.
By Rudin--Carleson theorem, we can find $H \in \cA(\bD)$ such that
\be
\|H\|_{\cH_\infty} = 1,
\qquad H(z) = - \f{V(z)}{|V(z)|},
\qquad H(\zeta_k) = V(\zeta_k), \quad k=1,\ldots,m.
\ee
Then, since $F=V-H \in \cA(\bD)$ satisfies $F(\zeta_1)=\cdots=F(\zeta_m)=0$, we have
\begin{align}
\mu_{\cA}(\ker(\cL_{\zeta}),\cV)
& \ge \f{\|F\|_{\cH_\infty}}{{\rm dist}_{\cH_\infty}(F,\cV)}
\ge \f{\|V-H\|_{\cH_\infty}}{\|H\|_{\cH_\infty}}
\ge |V(z)-H(z)| = |V(z)| \left(1+ \f{1}{|V(z)|} \right)\\
\nonumber
& = 1 + \eta.
\end{align}
This is the required inequality.
\epf

With Lemma \ref{LemObs} at our disposal, it becomes almost immediate to establish a result about optimal identification in $\cA(\bD)$ for acquisition functionals being evaluations at equispaced points.
There is a direct repercussion on optimal estimation in $\cA(\bD)$ --- the setting of Theorem \ref{ThmHInf} --- since
\be
\mu_{\cA}(\ker(\cL_{\zeta}),\cV)
= \sup_{\zeta_0 \in \bT} \mu_{\cA}(\ker(\cL_{\zeta}),\cV,e_{\zeta_0}).
\ee
The awaited result about optimal identification reads as follows.

\bprop
\label{PropMuHinfEqui}
Given points $\zeta_1,\ldots,\zeta_m \in \bT$ defined by
\be
\zeta_k =  \exp \left( i \f{2\pi }{m} (k-1) \right),
\qquad k=1,\ldots,m,
\ee
the indicator of compatibility of the acquisition process $\cL_{\zeta}$ and the approximability model relative to $\cP_n$ satisfies
\be 
\mu_{\cA}(\ker(\cL_{\zeta}),\cP_n)
= 2 + \kappa(m,n),
\qquad \mbox{with} \quad
\kappa(m,n) \asymp \ln \left( \f{m}{m-n+1} \right).
\ee
\eprop

\bpf
We invoke the result of \cite{RakShe}, which precisely says that
\be
\sup_{V \in \cP_n} \f{\|V\|_{\cH_\infty}}{\max\limits_{k=1,\ldots,m} |V(\zeta_k)|}
= 1 + \kappa(m,n),
\qquad \mbox{with} \quad
\kappa(m,n) \asymp \ln \left( \f{m}{m-n+1} \right).
\ee 
It remains to take Lemma \ref{LemObs} into account for the proof to be complete.
\epf

\begin{figure}[htbp]
	\center
	\subfigure[The indicator $\mu_{\cA}(\ker(\cL_{\zeta}),\cP_n,e_{\zeta_0})$ as a function of~$n$ for equispaced points and for random points.]{
		\includegraphics[width=0.45\textwidth]{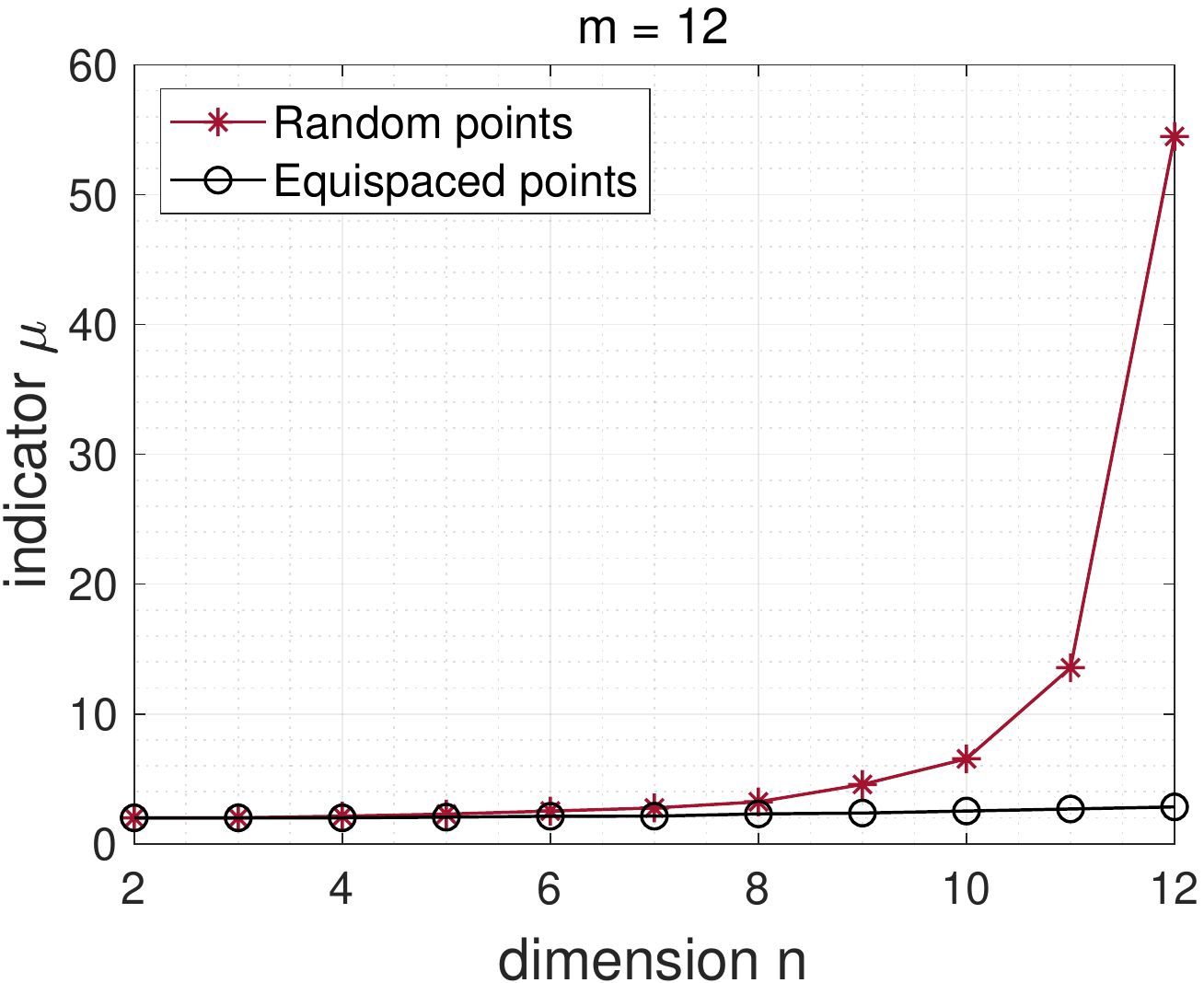}
	}
	\hspace{5mm}
	\subfigure[The estimation error $|F(\zeta_0)-A^{\rm opt}(\cL_\zeta(F))|$ as a \mbox{function} of $n$ for equispaced points and for random points.]{
		\includegraphics[width=0.45\textwidth]{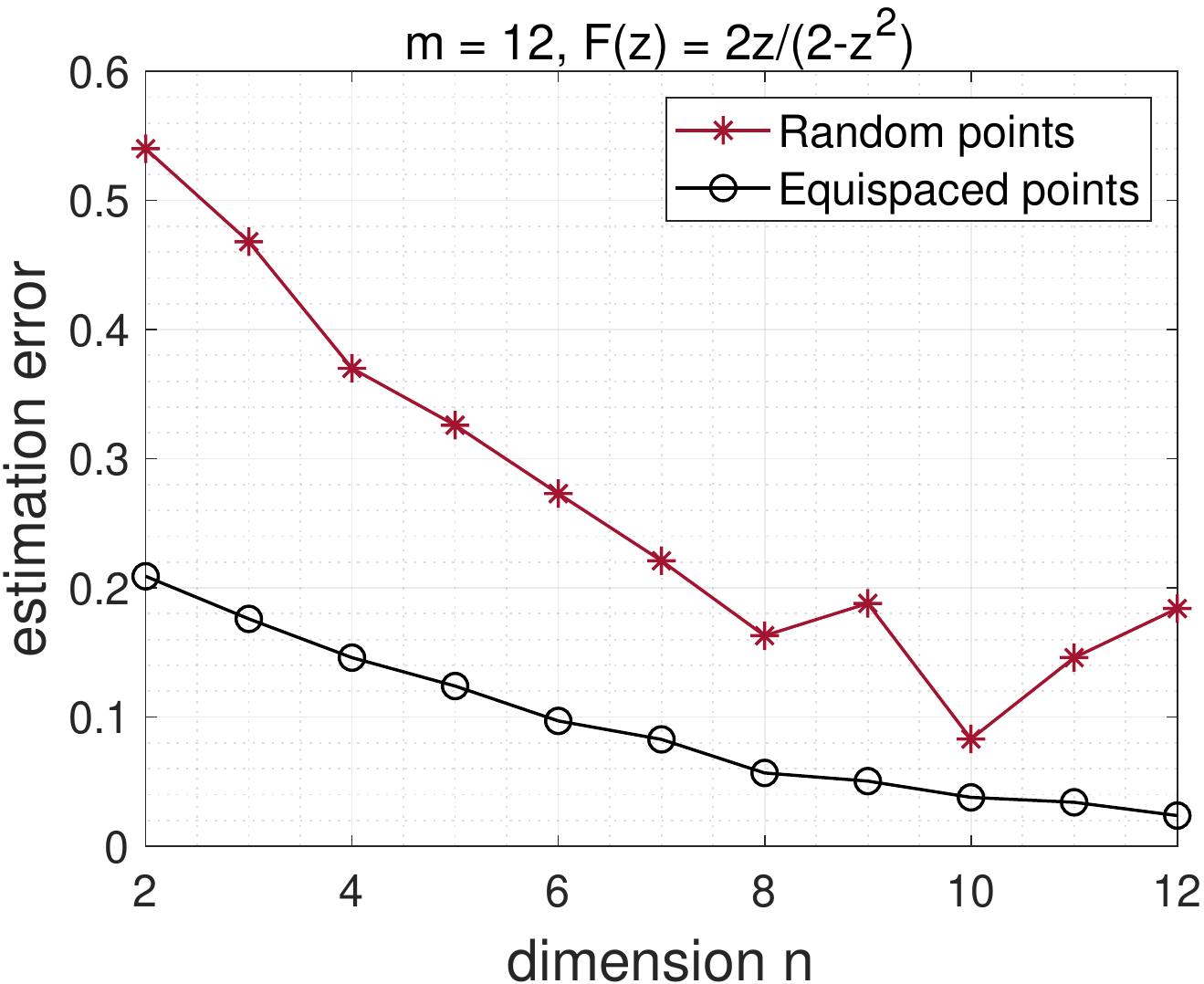}
	}
	\caption{Behavior of the indicator $\mu_{\cA}(\ker(\cL_{\zeta}),\cP_n,e_{\zeta_0})$ and of the estimation error \mbox{$|F(\zeta_0)-A^{\rm opt}(\cL_\zeta(F))|$} when the points $\zeta_0,\zeta_1,\ldots,\zeta_m$ are chosen on the torus.}
	\label{Fig2}
\end{figure}

To close this subsection, we present in Figure \ref{Fig2} 
a brief numerical illustration comparing the optimal estimation algorithms of Theorem \ref{ThmHInf} for equispaced points and random points on the torus.
The experiment is also included in the reproducible {\sc matlab} file.

\vspace{-2mm}
\section{Conclusion}
\label{SecConc}
\vspace{-5mm}

In this article, we formulated an optimal system identification problem by expressing the {\em a priori} information via approximability properties. To the best of our knowledge, this is the first work in this direction. We partially solved our problem by devising methods for the construction of optimal algorithms,
which turn out to be linear algorithms. Yet, our contribution is certainly an exploratory work that should be extended in several directions. We highlight a few possible directions below.

\vspace{-2mm}
{\em Measurement types.} 
We have mostly considered frequency responses obtained as point evaluations of transfer functions
and we have obtained sharp results when the evaluation points are equispaced. Our underlying method, however, is valid for any type of linear measurements, so it is worth studying its repercussions for other relevant measurement scenarios, including (i) impulse-response samples, (ii) convolutions of the impulse response with (pseudorandom) input signals, and (iii) pairs of input-output time series.

\vspace{-2mm}
{\em Measurements quality.}
The initial results presented here were based on the assumption that the observational data can be acquired with perfect accuracy.
In realistic situations, errors always occur in the acquisition process. It is possible to formulate an optimal-recovery problem taking this inaccuracy into account (in fact, such a problem can be transformed onto a standard optimal-recovery problem, at least formally --- see \cite{MicRiv} for details).
In the context of system identification, it would be valuable to design implementable algorithms that are optimal in this inaccurate setting.

\vspace{-2mm}
{\em MIMO systems.}
The focus in this preliminary work was put on single-input-single-output (SISO) systems. A step further would consist in treating multi-input-multi-output (MIMO) systems without simply applying our techniques to each input-output pair separately, 
as this becomes inefficient for large-scale systems.

\newpage

\section*{Appendix: Proofs of Essential Results}

In this section, we fully justify some statements appearing in the text but not yet established,
namely the 
\rev{relation between \eqref{Mod2}-\eqref{ModInf} and \eqref{KasIntersect},
	as well as}
the validity in the complex setting of results 
about optimal identification in Hilbert spaces~\cite{BCDDPW} and optimal estimation in Banach spaces~\cite{DFPW}.
We start with \rev{how \eqref{KasIntersect} connects to the descriptions \eqref{Mod2}-\eqref{ModInf} of the models put forward in \cite{HJN}.}

\bprop
With $\cX$ denoting either $\cH_2(\bD)$ or $\cA(\bD)$,
the following properties are equivalent:
\begin{align}
\label{Equiv1}
& \mbox{there exist } \rho>1 \mbox{ and } M>0 
\mbox{ such that } \|F(\rho \, \cdot)\|_{\cX} \le M;\\
\label{Equiv2}
& \mbox{there exist } \rho>1 \mbox{ and } M>0 
\mbox{ such that } 
{\rm dist}_{\cX}(F,\cP_n) \le M \rho^{-n}
\mbox{ for all } n \ge 0.
\end{align}
\eprop

\bpf
We write $F(z) = \sum_{n=0}^\infty f_n z^n$ throughout the proof.
We first establish the equivalence in the case $\cX = \cH_2(\bD)$.
Let us assume that \eqref{Equiv1} holds,
i.e., that $\sum_{n=0}^\infty |f_n|^2 \rho^{2n} \le M^2$ for some $\rho>1$ and $M>0$.
In particular, we have
$|f_n|^2 \le M^2 \rho^{-2n}$ for all $n \ge 0$.
It follows that, for all $n \ge 0$,
\be
{\rm dist}_{\cH_2}(F,\cP_n)^2 = \sum_{k=n}^\infty |f_k|^2 
\le M^2 \sum_{k=n}^\infty \rho^{-2k} = \f{M^2}{1-\rho^{-2}} \rho^{-2n},
\ee
hence \eqref{Equiv2} holds with a change in the constant $M$.
Conversely, let us assume that \eqref{Equiv2} holds,
i.e., that there are $\rho > 1$ and $M>0$ such that $\sum_{k=n}^\infty |f_k|^2 \le M^2 \rho^{-2n}$ for all $n \ge 0$.
In particular, we have $|f_n|^2 \le M^2 \rho^{-2n}$ for all $n \ge 0$.
Then, picking $\wt{\rho} \in (1,\rho)$, we derive that
\be
\sum_{n=0}^\infty |f_n|^2 {\wt{\rho}\,}^{2n}
\le M^2 \sum_{n=0}^\infty (\wt{\rho}/\rho)^{2n}
= \f{M^2}{1-(\wt{\rho}/\rho)^2},
\ee
hence \eqref{Equiv1} holds with a change in both $\rho$ and $M$.

We now establish the equivalence in the case $\cX = \cA(\bD)$.
Let us assume that \eqref{Equiv1} holds,
i.e., that  $\sup_{|z| = \rho} |F(z)| \le M$
for some $\rho > 1$ and $M>0$.
This implies that the Taylor coefficients of $F$ satisfy,
for any $k \ge 0$,
\be
|f_k| = \bigg|
\f{1}{2\pi i} \int_{|z|=\rho} \f{F(z)}{z^{k+1}} dz
\bigg|
\le \f{1}{2 \pi }\times \f{M}{\rho^{k+1}} \times 2 \pi \rho = M \rho^{-k}.
\ee
Considering $P \in \cP_n$ defined by $P(z):= \sum_{k=0}^{n-1} f_k z^k$,
we obtain
\be
{\rm dist}_{\cX}(F,\cP_n)
\le \|F-P\|_{\cH_\infty}
= \sup_{|z|=1} \bigg| \sum_{k=n}^\infty f_k z^k \bigg| 
\le \sum_{k=n}^\infty |f_k|
\le M \sum_{k=n}^\infty \rho^{-k} = \f{M}{1-\rho} \rho^{-n},
\ee
hence \eqref{Equiv2} holds with a change in the constant $M$.
Conversely, let us assume that \eqref{Equiv2} holds,
i.e., that there are $\rho > 1$ and $M>0$ such that
there exists, for each $n\ge 0$, 
a polynomial $P^{[n]} \in \cP_n$ with $\|F-P^{[n]}\|_{\cH_\infty} \le M \rho^{-n}$.
For all $n \ge 0$,
since the coefficients in $z^n$ of $F$ is the same as that of $F-P^{[n]}$, we have
\be
|f_n| = \bigg|
\f{1}{2 \pi i} \int_{|z|=1} \f{(F-P^{[n]})(z)}{z^{n+1}} dz
\bigg|
\le \f{1}{2 \pi} \times \|F-P^{[n]}\|_{\cH_\infty} \times 2 \pi 
\le M \rho^{-n}.
\ee
Then, picking $\wt{\rho} \in (1,\rho)$, we derive that
\be
\sup_{|z|=\wt{\rho}} |F(z)| \le \sum_{n=0}^\infty |f_n| {\wt{\rho}\,}^n
\le M \sum_{n=0}^\infty (\wt{\rho}/\rho)^n
= \f{M}{1-\wt{\rho}/\rho},
\ee
hence \eqref{Equiv1} holds with a change in both $\rho$ and $M$.
\epf

We now turn to the justification for the complex setting of the results from \cite{BCDDPW} about optimal identification in Hilbert spaces.
As in \cite[Theorem 2.8]{BCDDPW}, these results are easy consequences of the following statement.

\bthm
Let $\cV$ be a subspace of a Hilbert space $\cX$
and let $\ell_1,\ldots,\ell_m$ be linear functionals defined on $\cX$.
With a model set given by
\be
\cK = \{ f \in \cX: {\rm dist}_\cX(f,\cV) \le \eps \},
\ee
the performance of optimal identification from some $y \in \bC^m$ satisfies
\be
\inf_{A: \bC^m \to \cX} \sup_{f \in \cK \cap \cL^{-1}(y)}  \|f - A(y)\|_\cX
= \mu \left( \eps^2 - \min_{f \in \cL^{-1}(y)} \|f - P_\cV f\|_\cX^2  \right)^{1/2},
\ee
where the constant $\mu$ is defined by
\be
\mu = \sup_{u \in  \ker(\cL)} \f{\|u\|_\cX}{{\rm dist}_\cX(u,\cV)}.
\ee
\ethm

\bpf
Let $f^\star \in \cX$ be constructed from $y \in \bC^m$ via
$f^\star := \underset{f \in \cX}{\argmin} \|f - P_\cV f\|_\cX$
subject to 
$\cL(f) = y$.
We shall prove on the one hand that
\be
\label{Prosp1}
\sup_{f \in \cK \cap \cL^{-1}(y)} \|f - f^\star\|_\cX
\le \mu \left( \eps^2 - \|f^\star - P_\cV f^\star\|_\cX^2  \right)^{1/2}
\ee
and on the other hand that, for any $g \in \cX$,
\be
\label{Prosp2}
\sup_{f \in \cK \cap \cL^{-1}(y)} \|f - g\|_\cX 
\ge \mu \left( \eps^2 - \|f^\star - P_\cV f^\star\|_\cX^2  \right)^{1/2}.
\ee
Justification of \eqref{Prosp1}: Let us point out that $f^\star - P_\cV f^\star$ is orthogonal to both $\cV$ and $\ker(\cL)$.
To see this, given $v \in \cV$, $u \in \cU$, and $\theta \in [-\pi,\pi]$,
we notice that, as functions of $t \in \bR$,
the expressions
\begin{align}
\|f^\star - P_\cV f^\star + t e^{i \theta} v\|_\cX^2 
& = \|f^\star - P_\cV f^\star \|_\cX^2 + 2 t \Re( e^{- i \theta} \langle f^\star - P_\cV f^\star, v \rangle ) + \cO(t^2),\\
\|f^\star \hspace{-.5mm}+\hspace{-.5mm} t e^{i \theta} u - P_\cV( f^\star \hspace{-.5mm}+\hspace{-.5mm} t e^{ i \theta} u)\|_\cX^2 
& = \|f^\star - P_\cV f^\star \|_\cX^2 + 2 t \Re( e^{-i \theta} \langle f^\star - P_\cV f^\star, u - P_\cV u \rangle ) + \cO(t^2),
\end{align} 
are minimized at $t=0$.
Therefore $\Re( e^{- i \theta} \langle f^\star - P_\cV f^\star, v \rangle ) =0$
and $\Re( e^{- i \theta} \langle f^\star - P_\cV f^\star, u - P_\cV u \rangle ) = 0$ for all $\theta \in [-\pi,\pi]$.
This implies that $\langle f^\star - P_\cV f^\star, v \rangle  =0$ and  $\langle f^\star - P_\cV f^\star, u - P_\cV u \rangle =0$ 
for all $v \in \cV$ and $u \in \ker(\cL)$, hence our claim.
Now consider $f \in \cK \cap \cL^{-1}(y)$.
Since $\cL(f) = y = \cL(f^\star)$, we can write $f = f^\star + u $ for some $u \in \ker(\cL)$.
The fact that $f \in \cK$ then yields
\be
\eps^2 \ge \|f-P_\cV f\|_\cX^2 
= \|f^\star-P_\cV f^\star + u - P_\cV u\|_\cX^2
=  \|f^\star-P_\cV f^\star\|_\cX^2 + \| u - P_\cV u\|_\cX^2,
\ee
so that 
\be
{\rm dist}_\cX(u,\cV) = \| u - P_\cV u\|_\cX \le \left( \eps^2 -  \|f^\star - P_\cV f^\star\|_\cX^2  \right)^{1/2}.
\ee
It remains to take the definition of $\mu$ into account to obtain
\be
\| f - f^\star\|_\cX = 
\|u\|_\cX \le 
\mu \,  {\rm dist}_\cX(u,\cV) \le 
\mu \, \left( \eps^2 -  \|f^\star - P_\cV f^\star\|_\cX^2 \right)^{1/2}.
\ee
Justification of \eqref{Prosp2}: Let us select $u \in \ker(\cL)$ such that
\be 
\|u\|_\cX = \mu \, {\rm dist}_\cX (u,\cV)
\qquad \mbox{and} \qquad
\|f^\star - P_\cV f^\star\|_\cX^2 + \|u - P_\cV u\|_\cX^2 = \eps^2.
\ee
We now consider $f^\pm  := f^\star \pm u$.
It is clear that $f^\pm \in \cL^{-1}(y)$,
and we also have $f^\pm \in \cK$, since
\be
\|f^\pm - P_\cV f^\pm \|_\cX^2 
= \| (f^\star  - P_\cV f^\star)  \pm (u - P_\cV u ) \|_\cX^2 
=\|f^\star - P_\cV f^\star \|_\cX^2  + \|u - P_\cV u \|_\cX^2 
= \eps^2.
\ee
Then, for any $g \in \cX$,
\begin{align}
\sup_{f \in \cK \cap \cL^{-1}(y)} \|f - g\|_\cX & \ge \max_{\pm} \| f^\pm - g \|_\cX
\ge \f{1}{2} \left( \|f^+ - g\|_\cX + \|f^- - g\|_\cX \right) 
\ge \f{1}{2} \|f^+ - f^- \|_\cX \\
\nonumber
& = \|u\|_\cX = \mu \, {\rm dist}_\cX (u,\cV) = \mu \, \left( \eps^2 -  \|f^\star - P_\cV f^\star\|_\cX^2 \right)^{1/2}.
\end{align}
This completes the proof of the theorem.
\epf

Finally, we justify below that the result from \cite{DFPW} about optimal estimation in Babach spaces holds in the complex setting, too.

\bthm
Let $\cV$ be a subspace of a Banach space $\cX$, let $\ell_1,\ldots,\ell_m$ be linear functionals defined on $\cX$, and let $Q$ be another linear functional defined on $\cX$. 
With a model set given by
\be
\cK = \{ f \in \cX: {\rm dist}_\cX(f,\cV) \le \eps \},
\ee
the performance of optimal estimation of $Q$ satisfies
\be
\inf_{A: \bC^m \to \bC} \sup_{f \in \cK} \left| Q(f) - A(\cL(f)) \right|
= \mu \, \eps,
\ee
where the constant $\mu$ equals the minimum of the optimization problem
\be
\label{OptProg}
\minimize{a \in \bC^m} \bigg\| Q - \sum_{k=1}^m a_k \ell_k \bigg\|
_{\cX^*}
\qquad \mbox{subject to} \qquad
\sum_{k=1}^m a_k \ell_k(v) = Q(v) \quad \mbox{for all }v \in \cV.
\ee
\ethm

\bpf
Let $a^\star \in \bC^m$ be a minimizer of the optimization program \eqref{OptProg}
and let $\nu$ denote the value of the minimum.
Let us also consider
\be 
\mu = \sup_{u \in \cU} \f{|Q(u)|}{{\rm dist}_\cX(u,\cV)}.
\ee
We shall prove on the one hand that
\be
\label{Prosp3}
\sup_{f \in \cK} \bigg| Q(f) - \sum_{k=1}^m a^\star_k \ell_k(f) \bigg|
\le \nu \, \eps,
\ee
on the other hand that, for any $A: \bC^m \to \bC$,
\be
\label{Prosp4}
\sup_{f \in \cK} \left| Q(f) - A(\cL(f)) \right| \ge \mu \, \eps,
\ee
and we shall show as a last step that 
\be
\label{Prosp5}
\nu \le \mu.
\ee
Justification of \eqref{Prosp3}: 
Given $f \in \cK$, we select $v \in \cV$ such that $\|f-v\|_\cX = {\rm dist}_\cX(f,\cV)$.
The required inequality follows by noticing that
\begin{align}
\bigg| Q(f) - \sum_{k=1}^m a^\star_k \ell_k(f) \bigg|
& = \bigg| Q(f-v) - \sum_{k=1}^m a^\star_k \ell_k(f-v) \bigg|
\le \bigg\| Q- \sum_{k=1}^m a^\star_k \ell_k \bigg\|_{\cX^*} \|f-v\|_\cX\\
\nonumber
& = \nu \, {\rm dist}_\cX(f,\cV)
\le \nu \, \eps.
\end{align}
Justification of \eqref{Prosp4}: Let us select $u \in \ker(\cL)$ such that 
\be
|Q(u)| = \mu \, {\rm dist}_\cX(u,\cV)
\qquad \mbox{and} \qquad
\, {\rm dist}_\cX(u,\cV) = \eps.
\ee
Then, for any $A: \bC^m \to \bC$, we have
\begin{align}
\sup_{f \in \cK} |Q(f) - A(\cL(f))|
& \ge \max_\pm |Q(\pm u) - A(0)|
\ge \f{1}{2} \left( |Q(u)-A(0)| +  |Q(-u)-A(0)| \right)\\
\nonumber
& \ge \f{1}{2} |Q(u) - Q(-u)| = |Q(u)|
= \mu \, \eps.
\end{align}
Justification of \eqref{Prosp5}: 
We assume that $\ker(\cL) \cap \cV = \{ 0\}$,
otherwise $\mu = \infty$ and there is nothing to prove.
We consider a linear functional $\la$ defined on $\ker(\cL) \oplus \cV$ by
\be 
\la(u) = Q(u) \quad \mbox{for all }u \in \ker(\cL)
\qquad \mbox{and} \qquad
\la(v) = 0 \quad \mbox{for all }v \in \cV.
\ee
We then consider a Hahn--Banach extension $\wt{\la}$ of $\la$ defined on $\cX$.
Because $Q-\wt{\la}$ vanishes on $\ker(\cL)$,
we can write $Q-\wt{\la} = \sum_{k=1}^n \wt{a}_k \ell_k$ for some $\wt{a} \in \bC^m$,
and because $\wt{\la}$ vanishes on $\cV$, we have $\sum_{k=1}^m \wt{a}_k \ell_k(v) = Q(v)$ for all $v \in \cV$. 
We therefore derive that
\begin{align} 
\nu & \le \bigg\| Q - \sum_{k=1}^n \wt{a}_k \ell_k \bigg\|_{\cX^*}
= \big\|\wt{\la} \big\|_{\cX^*} = \| \la \|_{(\ker(\cL) \oplus \cV)^*}
= \sup_{\substack{u \in \ker(\cL)\\ v \in \cV}} \f{| \la(u-v) | }{\|u-v\|_\cX}\\
\nonumber 
&
= \sup_{\substack{u \in \ker(\cL)\\ v \in \cV}} \f{| Q(u) |}{\|u-v\|_\cX} = \sup_{u \in \ker(\cL)} \f{|Q(u)|}{{\rm dist}_\cX(u,\cV)} = \mu.
\end{align}
This concludes the proof of the theorem.
\epf

\end{document}